\documentclass[12pt,twoside]{article}
\usepackage{amsfonts}
\usepackage{amssymb,amsmath}
\begin{document}
\begin{center}
\textbf{The number of abc - equations $c^{n} = a+b$ satisfying
\\the strong abc - conjecture}
\bigskip
\\Constantin M. Petridi
\\ cpetridi@hotmail.com
\end{center}
\par
\vspace{2pt}
\begin{center}
\small{
\begin{tabular}{p{11cm}}
\textbf{Abstract.} We prove for the number $N(c,n)$ of equations\\
$c^{n}=a+b$, $(a,b,c,n)\,\in \mathbb{Z}^{4}_{+}$, $a<b$,
$(a,b)=1$, satisfying,\\ for any given $\varepsilon,\;
0<\varepsilon<1$, the strong abc - conjecture
$$c^{n}\,<\,R(c)^{\frac{\varepsilon}{1+\varepsilon}}R(a)^{\frac{1}{1+\varepsilon}}R(b)^{\frac{1}{1+\varepsilon}},$$
that $$\lim_{n\rightarrow
\infty}\,\frac{N(c,n)}{\frac{\phi(c^{n})}{2}}=1.$$
$$ $$
\bigskip
\end{tabular}
}
\end{center}
\par
\bigskip \textbf{1 Introduction}
\par\vspace{10pt}
In our paper [1] we proved that, for any given $0<\varepsilon<1$,
$$N(c)=\#\Bigl\{\;c<R(c)^{\frac{\varepsilon}{1+\varepsilon}}\,R(a_{i}b_{i})^{\frac{1}{1+\varepsilon}},
\;\,(a_{i},b_{i},c)\,\in \mathbb{Z}^{3}_{+},$$
$$c=a_{i}+b_{i},\;\,
a_{i}<b_{i},\,\;(a_{i},b_{i})=1;\,\;i=1,2,\ldots,
\frac{\phi(c)}{2}\;\Bigr\}$$ is equal to
$(1-\varepsilon)\frac{\phi(c)}{2}+O(\frac{\phi(c)}{2})$ for
$c\rightarrow \infty$. As usual $R(c)$ is the radical
$p_{1}p_{2}\ldots p_{\omega}$ of $c$,
$Q(c)=(p_{1}-1)(p_{2}-1)\ldots (p_{\omega}-1)$ and $\phi(c)$ is
the Euler function.
\par
\bigskip
In the present paper we restrict the equations to run only over
the powers $c^{n}$, $n=1,2,\ldots,\infty$. For the corresponding
number
$$N(c,n)=\#\Bigl\{\;c^{n}<R(c)^{\frac{\varepsilon}{1+\varepsilon}}\,R(a_{i}b_{i})^{\frac{1}{1+\varepsilon}},
\;\,(a_{i},b_{i},c,n)\,\in \mathbb{Z}^{4}_{+},$$
$$c^{n}=a_{i}+b_{i},\;\,
a_{i}<b_{i},\,\;(a_{i},b_{i})=1;\;\,i=1,2,\ldots,
\frac{\phi(c^{n})}{2}\;\Bigr\}$$ we prove the stronger result
$$\lim_{n\rightarrow
\infty}\,\frac{N(c,n)}{\frac{\phi(c^{n})}{2}}=1.$$
\par
\bigskip
The method used is based on our paper [2], where a lower bound and
an upper (trivial) bound are given for the geometric mean of the
radicals $R(ca_{i}b_{i})$ of the $\frac{\phi(c)}{2}$ equations
$c=a_{i}+b_{i},\;\,
0<a_{i}<b_{i}\,\;(a_{i},b_{i})=1;\;\,i=1,2,\ldots,
\frac{\phi(c)}{2}$, namely
$$\hspace{2.5cm}\kappa_{\varepsilon}\,\;R(c)^{1-\varepsilon}\,c^{2}\,<\,\Bigl\{\,\prod_{i=1}^{\frac{\phi(c)}{2}}\,
R(ca_{i}b_{i})\Bigr\}^{\frac{2}{\phi(c)}}\,<\,R(c)c^{2},\hspace{2.cm}(1)
$$ with $\kappa_{\varepsilon}$ an absolute positive constant, effectively
computable, depending only on $\varepsilon >0$.
\par
\bigskip\bigskip
\textbf{2 Limit Theorem}
\par
\bigskip
\textbf{Theorem.} For $(a_{i},b_{i},c,n) \in \mathbb{Z}^{4}_{+}$
and any given $\varepsilon$, $0<\varepsilon<1$, define $N(c,n)$ by
$$N(c,n)=\#\Bigl\{\;c^{n}<R(c)^{\frac{\varepsilon}{1+\varepsilon}}R(a_{i}b_{i})^{\frac{1}{1+\varepsilon}},\;\;
c^{n}=a_{i}+b_{i},\;a_{i}<b_{i},$$
$$(a_{i},b_{i})=1;\;\;i=1,2,\ldots,\frac{\phi(c^{n})}{2}\;\Bigr\}.$$
Then $$\lim_{n\rightarrow
\infty}\,\frac{N(c,n)}{\frac{\phi(c^{n})}{2}}=1.$$
\par
\bigskip
\textbf{Proof.} The relation
$c^{n}<R(c)^{\frac{\varepsilon}{1+\varepsilon}}R(a_{i}b_{i})^{\frac{1}{1+\varepsilon}}$
can also be written as
$$\hspace{2.3cm} R(c)^{1-\varepsilon}c^{n(1+\varepsilon)}<R(ca_{i}b_{i}),\;\;i=1,2,\ldots,\frac{\phi(c^{n})}{2}.\hspace{2.5cm}(2)$$
On the other hand since $c^{n}=a_{i}+b_{i}$, we have
$$\hspace{1.7cm}R(ca_{i}b_{i})=R(c)R(a_{i}b_{i})<R(c)c^{2n},\;\;i=1,2,\ldots,\frac{\phi(c^{n})}{2}\hspace{1.7cm}(3).$$
By definition therefore of $N(c,n)$, in the product
$$\Bigl(\prod_{i=1}^{\frac{\phi(c^{n})}{2}}R(ca_{i}b_{i}\Bigr)^{\frac{2}{\phi(c^{n})}},$$
$N(c,n)$ factors, in some order, are greater than
$R(c)^{1-\varepsilon}c^{n(1+\varepsilon)}$ because of (2), but are
smaller than $R(c)c^{2}$ because of (3). The remaining
$\frac{\phi(c^{n})}{2}-N(c,n)$ factors, according to same
definition of $N(c,n)$, are all smaller than
$R(c)^{1-\varepsilon}c^{n(1+\varepsilon)}$.
\par
\bigskip
In view of this and (1) it follows that
$$\kappa_{\varepsilon}R(c)^{1-\varepsilon}c^{2n}<
(R(c)c^{2n})^{N(c,n)\frac{2}{\phi(c^{n})}}\;(R(c)^{1-\varepsilon}c^{n(1+\varepsilon})^{1-N(c,n)\frac{2}{\phi(c^{n})}},$$
which, after making the necessary calculations in the exponents,
can be written as
$$\kappa_{\varepsilon}c^{n(1-\varepsilon)}<(R(c)^{\varepsilon}\;c^{n(1-\varepsilon)})^{N(c,n)\frac{2}{\phi(c^{n})}},$$
or by, taking logarithms
$$\log\kappa_{\varepsilon}+n(1-\varepsilon)\log c<(\varepsilon\log R(c)+n(1-\varepsilon)\log c){N(c,n)\frac{2}{\phi(c^{n})}}.$$
\par
\bigskip
Dividing by $\varepsilon\log R(c)+n(1-\varepsilon)\log c>0$ and
putting $c^{n}\frac{Q(c)}{R(c)}$ for $\phi(c^{n})$ we have,
$$\frac{\log\kappa_{\varepsilon}+n(1-\varepsilon)\log c}{\varepsilon\log R(c)+n(1-\varepsilon)\log c}
<N(c,n)\frac{2}{c^{n}\frac{Q(c)}{R(c)}}\leq 1,$$
\par
which, for $n\rightarrow\infty$, gives
$$\lim_{n\rightarrow\infty}\frac{N(c,n)}{\frac{\phi(c^{n})}{2}}=1,$$
or
$$N(c,n)=\frac{\phi(c^{n}\\)}{2}+o\Bigl(\frac{\phi(c^{n})}{2}\Bigr),\;\;n\rightarrow
\infty.$$ This completes the proof.
\par
\bigskip
\bigskip
\textbf{3 Note}
\par
\bigskip
The meaning of above Theorem is that almost all abc-equations
$c^{n}=a+b$ satisfy, for any given $0<\varepsilon<1$, the relation
$c^{n}<R(c)^{\frac{\varepsilon}{1+\varepsilon}}\,R(ab)^{\frac{1}{1+\varepsilon}}$.
\par
\bigskip
The method used to prove it, can be applied verbatim to the
abc-equation
$$c^{r}=a^{p}+b^{q},$$
\par
to show that the corresponding $N(c,r,p,q)$ is
$\frac{\phi(c^{r})}{2}+o\Bigl(\frac{\phi(c^{r})}{2}\Bigr),
\text{for } r \rightarrow \infty.$
\par
\bigskip
We wish to thank Peter Krikelis, Dep. of Mathematics,University of
Athens, for his help.
\par
\bigskip
\textbf{References}
\par
\bigskip
[1] Petridi, C.M., The number of equations satisfying the
abc-conjecture, arXiv:0904.1935 v1
\par
\bigskip
[2] Petridi, C.M., A strong "abc-conjecture" for certain
partitions $a+b$\\ of $c$, arXiv:math/0511224 v3

\end{document}